\documentclass[12pt]{amsart}

\usepackage{amsfonts}
\usepackage{amscd}
\usepackage{amssymb}

\allowdisplaybreaks

\usepackage{hyperref}

\date{\today}

\newtheorem{thm}{Theorem}[section]
\newtheorem{cor}[thm]{Corollary}
\newtheorem{lem}[thm]{Lemma}

\newtheorem{prop}[thm]{Proposition}
\theoremstyle{definition}

\theoremstyle{remark}
\newtheorem{rem}[thm]{Remark}

\numberwithin{equation}{section}

\newcommand{\R}{\mathbb R}
\newcommand{\m}{\bf{m}}
\newcommand{\Z}{\mathbb Z}
\newcommand{\T}{\mathbb T}
\newcommand{\Na}{\mathbb N}
\newcommand{\He}{\mathbb H}

\newcommand{\C}{{\mathbb C}}

\newcommand{\N}{\nabla }

\usepackage{color}

\title[Operators on Fock-Sobolev spaces]
{Fourier multipliers and pseudo-differential \\operators  
on Fock-Sobolev spaces }

\author[Sundaram Thangavelu ]{Sundaram Thangavelu }

\address[S. Thangavelu]{Department of Mathematics\\
Indian Institute of Science\\
560 012 Bangalore, India}

\email{veluma@iisc.ac.in}

\begin{document}

\maketitle

\vskip0.25in

\begin{abstract} Any bounded linear operator $ T $ on $ L^2(\R^n) $ gives rise to the operator $ S= B \circ T \circ B^\ast $ on the Fock space $ \mathcal{F}(\C^n) $ where $ B $ is the Bargmann transform. In this article we identify  those $ S $ which correspond to Fourier multipliers and pseudo-differential operators on $ L^2(\R^n)$ and study their boundedness on the Fock-Sobolev spaces $ \mathcal{F}^{s,2}(\C^n).$

 \end{abstract}


\section{Introduction} 

 Consider the Fock space $ \mathcal{F}(\C^n) $ consisting of entire functions $ F $ on $ \C^n$ which are square integrable with respect to the measure $ d \nu(z) = c_n e^{-\frac{1}{2}|z|^2} dz.$  In this article we are interested in certain bounded linear operators on the Fock space. The  well known Bargmann transform $ B $  is a unitary operator which takes  $ L^2(\R^n) $ onto  $ \mathcal{F}(\C^n) .$ Any  bounded linear operator $ T $ on $ L^2(\R^n) $ gives rise to the operator $ S = B \circ T \circ B^\ast $ on   $ \mathcal{F}(\C^n) .$ There is a  group of unitary operators $ \pi(z) $ indexed by $ \C^n $ acting on $ L^2(\R^n) $ by
$$ \pi(z)f(\xi) = e^{i(x\cdot \xi+\frac{1}{2} x \cdot y)} f(\xi+y),\,\, z =x+iy .$$
This family of unitary operators is related to the Schr\"odinger representation $ \pi_1 $ of the Heisenberg group $ \He^n,$  see Folland \cite{F}. As a consequence of the irreducibilty of  $ \pi_1 $ the only operators on $ L^2(\R^n) $ that commute with all $ \pi(z), z \in \C^n $ are constant multiple of the identity operator. On the other  hand there are families of operators which commute with $ \pi(x) $ or $ \pi(iy)$ for $x, y \in \R^n.$ For such an operator $ T $ the transferred operator $ S = B \circ T \circ B^\ast $ is expected to have certain invariance property.\\

 Observe that  $ \pi(ia) f(\xi) = f(\xi+a), a \in \R^n$ are translations and operators  $ T $ commuting with $ \pi(ia), a \in \R^n $ have explicit description. Indeed, all such operators are Fourier multipliers in the sense that  there exists $ m \in L^\infty(\R^n) $ such that 
$$ Tf(x) =: T_mf(x) = (2\pi)^{-n/2} \int_{\R^n} e^{i x \cdot \xi}\, m(\xi) \,\widehat{f}(\xi)\, d\xi.$$
Recall that the Bargmann transform (see Bargmann \cite{B} and Zhu \cite{Z1}) is defined by 
$$ Bf(z) = e^{-\frac{1}{4}z^2} \int_{\R^n} f(u) e^{-\frac{1}{2}|u|^2} e^{u \cdot z} du = e^{\frac{1}{4}z^2} \int_{\R^n} \widehat{f}(\xi) e^{-\frac{1}{2}|\xi|^2} e^{i \xi \cdot z} d\xi.$$
The translations $ \pi(ia) $ on $ L^2(\R^n) $ give rise to the action $ \rho(a) $ on $ \mathcal{F}(\C^n)$ defined by
$$ \rho(a) F(z) = F(z+a) e^{-\frac{1}{2} a \cdot z} e^{-\frac{1}{4} |a|^2} .$$ 
Thus we see that translation invariant operators on $ L^2(\R^n) $ correspond to operators on $ \mathcal{F}(\C^n) $ which commute with $ \rho(a), a\in \R^n.$ There is  another family of operators 
$ \pi(b)f(x) = e^{i b \cdot x} f(x),$  called the modulation operators,  which give rise to the operators $\rho(ib) $ acting on the Fock space by
$$ \rho(ib)F(z) = F(z+ib) e^{\frac{1}{2} ib \cdot z} e^{-\frac{1}{4} |b|^2} .$$
Operators on $ L^2(\R^n) $ which commute with $ \pi(b) $ are precisely (pointwise) multiplication operators. On the Fock space, they correspond to operators which commute with $ \rho(ib).$
More generally we have the following family of unitary operators indexed by $ \C^n$:
$$ \rho(w)F(z) = F(z+w) e^{-\frac{1}{2} \bar{w} \cdot z} e^{-\frac{1}{4} |w|^2} .$$
As these are related to an  irreducible unitary representations of the Heisenberg group realised on $\mathcal{F}(\C^n)$, the only operators on the Fock space that commute with $ \rho(w) $ for all $ w \in \C^n$ are constant multiples of the identity operator.\\

Consider   the Hilbert transform $ H $ on $ L^2(\R),$ which is a translation invariant operator given by the multiplier $ m(\xi) = -i \, sign \,\xi.$  In  \cite{Z} Zhu  calculated the transferred operator and showed that
$$ (B \circ H \circ B^\ast) F(z) = \int_{\C} F(w)\, \varphi(z-\bar{w}) \,e^{\frac{1}{2} z \bar{w}} \, d\nu(w)$$
where $ \varphi$ is explicitly given in terms of the error function $ A.$ Motivated by this example, Zhu suggested considering operators of the form
$$ S_\varphi F(z) =  \int_{\C^n} F(w)\, \varphi(z-\bar{w}) \,e^{\frac{1}{2} z \cdot \bar{w}} \, d\nu(w)$$
where $ \varphi \in \mathcal{F}(\C^n) $ and raised the question of characterising those $ \varphi $ for which $ S_\varphi $ are bounded on $ \mathcal{F}(\C^n) .$ Using the above mentioned connection between Fourier multipliers, the authors in \cite{CLSWY} proved that $ S_\varphi $ is bounded on $ \mathcal{F}(\C^n) $ if and only if 
$$ \varphi(z) = \int_{\R^n} m(\xi) e^{-(\xi-\frac{i}{2}z)^2} d\xi$$ fr some $ m \in L^\infty(\R^n).$ A shorter proof of this can be given by calculating the action of $ B \circ T_m \circ B^\ast $ to the reproducing kernel $ k(z, \bar{w}) = e^{\frac{1}{2} z\cdot \bar{w}}.$ 
In a similar fashion the operators on $ \mathcal{F}(\C^n) $ commuting with $ \rho(ib) $ are of the form
$$ \widetilde{S}_\varphi F(z) = \int_{\C^n} F(w) e^{\frac{1}{2}z \cdot \bar{w}} \, \varphi(z+\bar{w}) \, d\nu.$$
These operators correspond to pointwise multipliers on $ L^2(\R^n)$ and have already been studied by Bais and Dogga \cite{BD}. We can also describe $ \varphi $ in terms of $ m $ which is defined by the relation $ (B^\ast \circ \widetilde{S}_\varphi \circ B) f = m f.$\\

Though the operators $ S_\varphi $ and $ \widetilde{S}_\varphi $ look deceptively similar, their behaviours are very different. In fact there is an uncertainty principle on the boundedness of these operators.

\begin{thm} For $ \varphi \in \mathcal{F}(\C^n), $  the operator $ S_\varphi $ is bounded on the Fock space if and only if 
$ \varphi(z) = \int_{\R^n} m(\xi) e^{-(\xi-\frac{i}{2}z)^2} d\xi$ 
for some $ m \in L^\infty(\R^n)$ whereas  $ \widetilde{S}_\varphi $ is bounded if and only if 
$ \varphi(z) = \int_{\R^n} m(\xi) e^{-(\xi-\frac{1}{2}z)^2} d\xi $ 
for some other $ m \in L^\infty(\R^n).$  Moreover, for any non constant $ \varphi$ either $ S_\varphi $ or $ \widetilde{S}_\varphi $ is unbounded.
\end{thm}

The description of $ \varphi $ in the above theorem suggests that we consider the following transform. For $ f \in L^2(\R^n,d\gamma)$ where  $ d\gamma(x) = e^{-|x|^2} dx ,$ define
$$ Gf(z) = \int_{\R^n} f(\xi) e^{-(\xi-\frac{i}{2}z)^2} d\xi = e^{\frac{1}{4} z^2} \, \int_{\R^n} f(\xi) e^{i z \xi} \, d\gamma(\xi).$$
We can express $ G $ in terms of the Bargmann transform and for that reason it is known as the Gauss-Bargmann transform.  Define the unitary operator $ U $ on the Fock space by 
$ UF(z) = F(-iz).$ Then we see that $ Gf = (U^\ast \circ B \circ E) f $ where $  Ef(\xi) = f(\xi) e^{-\frac{1}{2}|\xi|^2}.$ As $ E $  is a unitary operator from  $L^2(\R^n,d\gamma)$ onto 
$ L^2(\R^n)$ it follows that $ G :L^2(\R^n,d\gamma) \rightarrow \mathcal{F}(\C^n) $ is unitary.   Thus we see that $ S_\varphi $ is bounded if and only if $ \varphi = Gm $ for some $ m \in L^\infty(\R^n) \subset L^2(\R^n,d\gamma).$ In other words, $ \varphi $ has to belong to the subspace $ G L^\infty(\R^n) \subset \mathcal{F}(\C^n).$ Working with $ G $ has another advantage: we have the following interesting description of the operators $ S_\varphi.$\\

\begin{thm} A bounded linear operator $ S $ on $ \mathcal{F}(\C^n) $ commutes with all $ \rho(a), a \in \R^n$  if and only if there exists an $ m \in L^\infty(\R^n) $ such that
$$ SF(z) = e^{\frac{1}{4}z^2} \, \int_{\R^n} m(\xi)\, G^\ast F(\xi)\, e^{i z \cdot \xi}\, d\gamma(\xi).$$
\end{thm}

\begin{rem} The above representation suggests that we use the notation $ S_\varphi F= \varphi \ast F .$  Then $ G^\ast $ takes this convolution into products: $ G^\ast(\varphi \ast F ) = (G^\ast \varphi)\, (G^\ast F) $ as $ \varphi = Gm.$ There is a similar representation for the operators $ \widetilde{S}_\varphi.$
\end{rem}

\begin{thm} A bounded linear operator $ \widetilde{S} $ on $ \mathcal{F}(\C^n) $ commutes with all $ \rho(ib), b \in \R^n$  if and only if there exists an $ m \in L^\infty(\R^n) $ such that
$$ \widetilde{S}F(z) = e^{\frac{1}{4}z^2} \, \int_{\R^n} m(\xi)\, \, B^\ast F(\xi)\, e^{ z \cdot \xi}\, d\xi.$$
\end{thm}

\begin{rem}We note that in the above theorem the multiplier $ m $ corresponding to $ \widetilde{S}_\varphi $ is given by $ m(\xi) e^{-\frac{1}{2}|\xi|^2} = B^\ast \varphi.$
\end{rem}

In a recent work \cite{CHLS} the boundedness of operators of the form $ S_\varphi$ have been investigated on fractional order Fock-Sobolev spaces $ \mathcal{F}^{s,2}(\C^n), s \geq 0 $ which are the images under Bargmann transform of the Hermite-Sobolev spaces $ W_H^{s,2}(\R^n).$ Let us recall that with $ H = -\Delta+|x|^2 $ denoting the Hermite operator, the space $ W_H^{s,2}(\R^n)$ is the subspace of $ L^2(\R^n) $ consisting of those $ f $ for which $ H^{s/2}f \in L^2(\R^n).$ Here the fractional powers $ H^{s/2} $ are defined using the explicit spectral decomposition of the Hermite operator.  It is also known from the work of Cho and Park \cite{CP} that $ \mathcal{F}^{s,2}(\C^n)$ are also weighted Bergman spaces. As Hermite functions are mapped onto the monomials $ z^\alpha $ the spaces $ \mathcal{F}^{s,2}(\C^n)$ are natural to look at. Apart from being translation invariant, the spaces $ W_H^{s,2}(\R^n)$ are also invariant under the Fourier transform and hence every Fourier multiplier on these spaces is also a pointwise multiplier. Using these properties, the authors in \cite{CHLS} proved that $ S_\varphi $ is bounded on $ \mathcal{F}^{s,2}(\C^n)$ if and only if $ \varphi $ has the same representation as  in Theorem 1.1 where now $ m $ defines a Fourier multiplier on $ W_H^{s,2}(\R^n).$  When $ s = 0 $ it is known that $ S_\varphi $ is bounded on $ \mathcal{F}(\C^n)$ if and only the associated $ m $ belongs to $ L^\infty(\R^n).$ For the fractional case, there is no such characterisation is known.\\

The image of $ L^2(\R^n, d\gamma) $ under $ G $ is the same as the image of $ L^2(\R^n) $ under $ B $ which is $ \mathcal{F}(\C^n).$  For any $ s \geq 0,$ 
$ W_H^{s,2}(\R^n) $  is the image of $ L^2(\R^n) $ under $ H^{-s/2}.$ Similarly, we can define the Gauss-Sobolev spaces $ W_G^{s,2}(\R^n) $ as the image of  $ L^2(\R^n, d\gamma) $ under $ L^{-s/2} $ where $ L = -\Delta + x \cdot \N$ is the Ornstein-Uhlenbeck operator. It turns out that the image of $ W_G^{s,2}(\R^n) $ under $G$ coincides with the image of $ W_H^{s,2}(\R^n) $ under $ B ,$ namely with  $ \mathcal{F}^{s,2}(\C^n) .$
The boundedness of the operators $ S_\varphi$ on $ \mathcal{F}^{s,2}(\C^n) $ has been studied by Wick and Wu in \cite{WW} when $ s = k $ is a non negative integer and the general case 
by Cao et al in \cite{CHLS}. We  can simplify their proofs and improve the results by obtaining the following characterization. Let $ MW_H^{s,2}(\R^n) $ be the set of all Fourier multipliers of the Hermite-Sobolev space $ W_H^{s,2}(\R^n).$ We denote by  $ M\mathcal{F}^{s,2}(\C^n) $  the image of $ MW_H^{s,2}(\R^n) $ under $ G $ which is a proper subspace of $ \mathcal{F}^{s,2}(\C^n).$ The results of \cite{WW} and \cite{CHLS} can be restated as follows. 

\begin{thm} The operator $ S_\varphi $ is bounded on  $ \mathcal{F}^{s,2}(\C^n) $ if and only if $ \varphi \in  M\mathcal{F}^{s,2}(\C^n).$  In other words,  $ \varphi = Gm $ for some $ m \in  MW_H^{s,2}(\R^n) .$ \\
\end{thm}


Since we do not have a characterisation of $ MW_H^{s,2}(\R^n), $ it is preferable to find sufficient conditions on $ m $ which are easy to check. As a partial characterisation for the case  $ s =k$, we can show that $ m \in MW_H^{k,2}(\R^n)$ if and only if $ \partial^\alpha m \circ H^{-|\alpha|/2} $ is bounded on $ L^2(\R^n) $ for all $ |\alpha| \leq k.$ \\

A useful subclass of multipliers is provided by  $ L_k^\infty(\R^n), $  the space of functions which satisfy the 
estimates $ |\partial^\alpha m(x) |\leq C (1+|x|)^{|\alpha|} $ for all $ |\alpha| \leq k.$ It can be shown that $ L_k^\infty(\R^n) $ is contained in $ MW_H^{k,2}(\R^n) .$ We let 
$ \mathcal{F}_0^{k,2}(\C^n) $ stand for the image of $L_k^\infty(\R^n) $ under the Gauss-Bargmann transform. When $ k = 0 $ we simply write $ \mathcal{F}_0(\C^n) $ instead of $ \mathcal{F}_0^{0,2}(\C^n) .$

\begin{thm} Let $ k $ be a non-negative integer. If we assume that $ \varphi \in \mathcal{F}_0^{k,2}(\C^n) \subset \mathcal{F}^{k,2}(\C^n) ,$  or equivalently $ \varphi = Gm, m \in L_k^\infty(\R^n),$ then $ S_\varphi $ is bounded on $ \mathcal{F}^{k,2}(\C^n) .$ \\
\end{thm}

{\bf Example.} An example of $ \varphi  \in \mathcal{F}_0^{k,2}(\C^n) $ is provided by the function $ \varphi_t(z) = c_t \, e^{-\frac{1}{4} \frac{it}{1+it} z^2} $ which corresponds to the multiplier $ m_t(\xi) = e^{-it |\xi|^2}, \, t \in \R.$ The Fourier multiplier $ T_{m_t} $ is the Schr\"odinger group $ e^{it \Delta} $ generated by the Laplacian $ \Delta$ on $\R^n.$ Thus the operator $ B \circ e^{it\Delta} \circ B^\ast $ is bounded on every $ \mathcal{F}^{k,2}(\C^n).$ The function $ u(z,t) = e^{-\frac{1}{4} z^2} S_{\varphi_t}F(z) $ solves the Schr\"odinger equation associated to the complexified Laplacian $ \Delta_\C = \sum_{j=1}^n \frac{\partial^2}{\partial z_j^2} $ on $ \C^n.$\\

\begin{rem} As the Hermite-Sobolev spaces are invariant under the Fourier transform, it follows that $ \partial^\alpha m \circ H^{-|\alpha|/2} $ is bounded on $ L^2(\R^n) $ if and only if $ T_{\partial^\alpha m} \circ H^{-|\alpha|/2} $ is bounded on $ L^2(\R^n) .$ As the image of the Hermite operator $ H $ under the Bargmann transform is the number operator $ \mathcal{R} = 2 \sum_{j=1}^n z_j \frac{\partial}{\partial z_j} +n $ on the Fock space, it follows that the above condition on $ \partial^\alpha m$ is equivalent to the boundedness of $ S_{\partial^\alpha \varphi} \circ \mathcal{R}^{-|\alpha|/2}$ on $ \mathcal{F}(\C^n)$ where $ \varphi = Gm.$
 It would be interesting to find sharp conditions on $ \varphi $ which allows us to conclude that $ S_\varphi $ is bounded on the Fock-Sobolev spaces.\\
\end{rem}

We now consider bounded linear operators  on $ \mathcal{F}(\C^n) $ which are not necessarily of the form $ S_\varphi.$ Note that any such operator $ S $ on the Fock space   is an integral operator of the form 
$$ SF(z) = \int_{\C^n} F(w) K(z,\bar{w}) d\nu(w) $$ for a suitable  kernel $ K(z,w) = K_S(z,w) $ which is a holomorphic function on $ \C^{2n}.$ Indeed, as we have the reproducing formula
$$ F(z) = \int_{\C^n} F(w) g(\cdot, \bar{w}) d\nu(w),\,\,  g(z,w) = g_w(z) =: e^{\frac{1}{2} z\cdot w} $$
it follows that $ K_S(z,w) = Sg(\cdot,w)(z).$ This suggests that one considers  integral operators  with kernel $ K(z,w) $ and ask for sufficient conditions under which they are bounded on $ \mathcal{F}(\C^n).$ Slightly changing the notation and writing $ K(z,w) = k(z,z-w) e^{\frac{1}{2} z \cdot w} $ we may also write them in  the form
$$ S_KF(z) = \int_{\C^n} F(w) k(z,z-\bar{w}) \, e^{\frac{1}{2} z\cdot \bar{w}} \, d\nu(w) $$
and call them pseudo-differential operators on the Fock space for the following reason.\\

Recall that when  $ m = \sum_{|\alpha| \leq k} a_\alpha  (-i)^{|\alpha|}\, \xi^\alpha $ is a polynomial, the Fourier multiplier operator $ T_m $ reduces to the differential operator $ m(D) = \sum_{|\alpha| \leq k} a_\alpha  \, \partial^\alpha. $  Since we can also write the polynomial $ m $ in terms of the Hermite polynomilas as $ m(\xi) = \sum_{|\alpha| \leq k} b_\alpha  (-i)^{|\alpha|}\, H_\alpha(\xi) $ the corresponding operator $S_\varphi $ on the Fock space takes the form
$$ S_\varphi F(z) = e^{\frac{1}{2}z^2} \int_{\R^n} p(\xi) \, G^\ast F(\xi)\, e^{i \xi \cdot z} \, d\gamma = p\ast F(z) $$
where $ p(z) = \sum_{|\alpha| \leq k} b_\alpha  z^\alpha.$
Thus we see that partial differential operators on $ L^2(\R^n) $ give  rise to convolution operators with polynomials. Note that the operators $ S_\varphi $ correspond to kernels $ k(z,w) $ which are independent of the first variable. It is therefore natural to consider operators on the Fock space which correspond to pseudo-differential operators  $ a(x,D) $ defined on $ L^2(\R^n) $ by
$$ a(x,D)f(x) = (2\pi)^{-n/2} \int_{\R^n} a(x,\xi)\, \widehat{f}(\xi)\, e^{i x \cdot \xi}\, d\xi.$$
By doing a formal calculation, these operators can also be written in the form  where $ \mathcal{F}_2 $ is the Fourier transform in the second set of variables.
$$ a(x,D)f(x) = (2\pi)^{-n/2} \, \int_{\R^n} \mathcal{F}_2^{-1}a(x,x-y)\, f(y) \, dy $$ 
The analogy between $ S_K $ and $ a(x,D) $ is now clear.\\

However, the operators $ S_K $ are more general in the sense that $ B^\ast \circ S_K \circ B $ need not be a pseudo-differential operator on $ L^2(\R^n).$ 
We would like to investigate the operators $ S_K $ which actually correspond to pseudo-differential operators. In other words, we would like to identify the kernels $ K $ so that $ S_K =: B^\ast \circ a(x,D) \circ B $  for a symbol $ a(x,\xi).$  \\

 \begin{rem} In view of the representation for $ S_\varphi$ established in Theorem 1.2,  there is yet another possibility of defining pseudo-differential operators on Fock spaces. Given a symbol $b(z,\xi) $ which is entire in the first variable, we can consider
$$ b(z,D)F(z) = e^{\frac{1}{2}z^2} \int_{\R^n}  b(z, \xi) \, G^\ast F(\xi)\, e^{i \xi \cdot z} \, d\gamma .$$
It would be interesting to find out the connection between these two families of pseudo-differential operators.
\end{rem}

 The pseudo-differential operators $ a(x,D) $ defined above correspond to Kohn-Nirenberg class. There is another equivalent definition which correspond to the  Weyl calculus, given as  Weyl transforms of tempered distributions. 
The Weyl transform of a function $ \sigma \in L^1(\C^n) $ or a finite Borel measure $ \mu $ on $ \C^n $ we define
$$ W(\sigma)f(\xi) = \int_{\C^n} \sigma(z) \pi(z)f(\xi) \, dz.$$ 
It is known that $ W(\sigma) $ can be written as $ \widehat{a}(x,D) $ for a suitable symbol $ a.$ The transferred operators $ B \circ W(\varphi) \circ B^\ast $ are easy to calculate.  In fact, a simple calculation shows that the action of $ \rho(w) = B \circ \pi(w) B^\ast $ on $ F \in \mathcal{F}(\C^n) $ is given by 
$$ \rho(w)F(z) = F(z+w) e^{-\frac{1}{2} z \cdot \bar{w}} e^{-\frac{1}{4} |w|^2} $$
and hence we see that $(B \circ W(\varphi) \circ B^\ast ) F(z) $ is given by the integral
$$ \int_{\C^n} \varphi(w) \rho(w)F(z) dw = \int_{\C^n}  \varphi(w) F(z+w) e^{-\frac{1}{2} z \cdot \bar{w}} e^{\frac{1}{4} |w|^2} d\nu(w).$$
Thus the boundedness of $ W(\varphi) $ on $ L^2(\R^n)$ is equivalent to the boundedness of the operator appearing on the right hand side of the above equation. However, it is more informative to look at the integral representation of such operators with holomorphic kernels. \\

We now specialise to  integral operators on the Fock space with kernels 
$ K(z,w) =  B\psi(z,w) $ for some  $ \psi $ on $\R^{2n} .$ Given $ \sigma \in L^1(\C^n) $ let us write $ \sigma(x,y) = \sigma(x+iy) $ and define $$ \psi(u,v) = (\mathcal{F}_1^{-1}\sigma)( \frac{1}{2}(u+v), u-v) $$ where $ \mathcal{F}_1 $ is the Fourier transform in the first set of variables. Note that we can recover $\sigma $ from  $ \psi $ from the formula
\begin{equation}\label{sigma}  \sigma(z) = (2\pi)^{-n/2}\, \int_{\R^n} e^{-i x\cdot \xi} \psi(\xi+\frac{1}{2}y, \xi-\frac{1}{2}y) d\xi.
\end{equation}

\begin{thm} Given  a function $ \psi $ on $ \R^{2n},$  let $ \sigma $ be defined as above. Then the operator
$$ S_K F(z) = \int_{\C^n} F(w) \, B\psi(z,\bar{w}) \, d\nu(w) $$
corresponds to the  pseudo-differential operator  $ W(\sigma)$ on $ L^2(\R^n).$  Consequently, $ S_K $ is bounded on $ \mathcal{F}^{s,2}(\C^n)$ if and only if $ W(\sigma) $ is bounded on the Hermite -Sobolev space $ W_H^{s,2}(\R^n).$ \\
\end{thm}

 When $ \psi \in L^2(\R^{2n}),$ the function $ \sigma \in L^2(\R^{2n}) $ and hence both $ W(\sigma) $ and $S_K $ are Hilbert-Schmidt operators. Note that in this case $ B\psi  \in \mathcal{F}(\C^{2n}).$ The boundedness of $ W(\sigma) $ on $ W_H^{s,2}(\R^n) $ is equivalent to the boundedness of $ H^{s/2} \circ W(\sigma)\circ H^{-s/2} $ on $ L^2(\R^n).$ As each one of these operators is a pseudo-differential operator, we can appeal to some general theorems to find a sufficient condition on $ \sigma$ so that $ S_K $ is bounded on $ \mathcal{F}^{s,2}(\C^n),$ see Remark 4.1.\\
 
 A particularly interesting class of operators emerge when we choose $ \sigma $ to be a radial function/measure on $ \R^{2n}.$ In this case the operator $ W(\sigma) $ is diagonalised by the Hermite functions and we can explicitly compute the kernel $ K(z,w) $ associated to the transferred operator, see Proposition 4.2. The kernel turned out to be of the form $ k(\sqrt{z \cdot w}) $ for an even entire function $ k $ of a single complex variable.\\

As the Bargmann transform is intimately connected to the theory of Hermite expansions, we studied  translation invariant operators on $ W_H^{s,2}(\R^n) $ and their images under the Bargmann transform.  These spaces, as well their images, namely the fractional Fock-Sobolev spaces $ \mathcal{F}^{s,2}(\C^n)$ are invariant under the Fourier transform. The classical Sobolev spaces $ W^{s,2}(\R^n), s \geq 0 $ form another class of translation invariant subspaces  of $ L^2(\R^n) $ and hence suitable for the study of Fourier multipliers. It is therefore natural to study the images of these spaces under the Bargmann transform, denoted by $ \mathcal{F}_{1,0}^{s,2}(\C^n)$, and study the boundedness of the  operators $ S_\varphi $ corresponding to the Fourier multipliers $ T_m.$ As $ W^{s,2}(\R^n)$ are not invariant under  the Fourier transform, unlike the case of $ \mathcal{F}^{s,2}(\C^n),$ we cannot expect  $ \mathcal{F}_{1,0}^{s,2}(\C^n)$ to be described by  radial weight functions. However, it is possible to describe them  differently using derivatives. Some results on the boundedness of $ S_\varphi$ on these spaces are proved in Section 3.\\

The plan of the paper is as follows. In Section 2 we set up the notations and recall relevant facts regarding Bargmann transform, Hermite, Gauss and Fock-Sobolev spaces. In Section 3, we prove the results on operators $ S_\varphi$ whereas in Section 4 we treat the operators $ S_K $ which correspond to pseudo-differential operators on $ L^2(\R^n).$\\

\section{Preliminaries}
 
 \subsection{Hermite functions and the Bargmann transform} By Taylor expanding the entire function $  g_x(w) = e^{-w^2+2xw} $ into power series we get
 $$ e^{-w^2+2xw} = \sum_{k=0}^\infty \frac{1}{k!} H_k(x)\, w^k .$$
 Since  $ g_x(w) = e^{x^2} e^{-(x-w)^2} $ it is clear that $ H_k(x) = \frac{d^k}{dx^k}g_x(0) $ is a polynomial of degree $ k.$ These are called the Hermite polynomials. The Hermite polynomials $ H_\alpha(x), \alpha \in \Na^n $ on $ \R^n $ are defined by taking the products of the one dimensional Hermite functions. The normalised Hermite functions are then defined by
 $$ \Phi_\alpha(x) = ( \pi^{n/2} 2^{|\alpha|} \alpha !)^{-1/2} H_\alpha(x)\, e^{-\frac{1}{2}|x|^2} .$$
 It is known that this family $ \Phi_\alpha, \alpha \in \Na^n $ forms an orthonormal basis for $ L^2(\R^n).$  By defining $ \zeta_\alpha(w) =  ( \pi^{n/2} 2^{|\alpha|} \alpha !)^{-1/2} w^\alpha$ 
 the generating function identity for the Hermite functions takes the form 
  $$   \sum_{\alpha \in \Na^n} \Phi_\alpha(x)\, \zeta_\alpha(w)  = \pi^{-n/2} e^{\frac{1}{4}w^2} e^{-\frac{1}{2}(x-w)^2}.$$
 We  use the above identity in defining the Bargmann transform and deduce  many of its properties easily. For the properties of the Bargmann transform we  need we refer to the paper \cite{B} of Bargmann; see also the book \cite{Z1} by Zhu. For the results on Hermite functions we refer to Szego \cite{Sz} and the monograph \cite{T1}.\\
 
 For $ f \in L^2(\R^n) $ we define $ Bf $  by integrating $f $ against the kernel on the right hand side:
 $$ Bf(w) = \pi^{-n/2} e^{\frac{1}{4} w^2} \int_{\R^n} f(x) \,e^{-\frac{1}{2}(x-w)^2} dx =  \pi^{-n/2} e^{-\frac{1}{4} w^2} \int_{\R^n} f(x) \,e^{-\frac{1}{2}x^2}\, e^{x\cdot w} dx.$$
 In view of the above identity, it is clear that $ B\Phi_\alpha = \zeta_\alpha $  and we have the expansion 
 $$ Bf(w) =  \sum_{\alpha \in \Na^n} (f, \Phi_\alpha)\, \zeta_\alpha(w). $$ 
 If we let $ d\nu(w) = c_n e^{-\frac{1}{2}|w|^2} $ then the functions $ \zeta_\alpha $ form an orthonormal system in $ L^2(\C^n, d\nu).$ By the Plancherel theorem for the Hermite expansions we deduce that
 $$ \int_{\C^n} |Bf(w)|^2\, d\nu(w) =  \sum_{\alpha \in \Na^n} |(f, \Phi_\alpha)^2\, = \int_{\R^n} |f(x)|^2\, dx.$$
 Thus for $ f \in L^2(\R^n) $ its Bargmann transform $ Bf $ is an entire function which is square integrable with respect to $ d\nu.$ It can be shown that any entire function $ F $ on $ C^n $ which is in $ L^2(\C^n, d\nu) $ is of the form $ Bf $ for a unique $ f \in L^2(\R^n).$ Therefore, if we define the Fock space $ \mathcal{F}(\C^n) $ as the space of $ L^2(\C^n, d\nu) $ consisting of entire functions, then $ B: L^2(\R^n) \rightarrow \mathcal{F}(\C^n) $ is unitary. The inverse of $ B $ is given by its adjoint $ B^\ast.$\\

\subsection{Gauss-Bargmann transform} Along with the Bargmann transform $ B $ which takes $ L^2(\R^n) $ onto $ \mathcal{F}(\C^n) $ we also consider the Gauss-Bargmann transform which is more suitable for the study of the operators $ S_\varphi.$ Let $ UF(z) = F(-iz) $ which is a unitary operator on $ \mathcal{F}(\C^n) $ which intertwines the Fourier transform: $ B^\ast \circ U \circ B = \mathcal{F}.$  Let $ d\gamma(x) = e^{-|x|^2} dx $ and consider the operator $ G $ defined on $ L^2(\R^n, d\gamma) $ by 
$$ Gf(z) = e^{\frac{1}{4}z^2} \int_{\R^n} f(x) e^{i x\cdot z} \, d\gamma(x) = \int_{\R^n} f(x) e^{-(x-\frac{i}{2}z)^2} \, dx.$$
Note that $ G $ is related to $B $ via  $ Gf(z) = (U\circ B)g(z) $ where $ g(x) = f(x) e^{-\frac{1}{2}|x|^2}.$ We can read out the properties of $ G $ from that of $ B.$ Thus $ G: L^2(\R^n, d\gamma) 
\rightarrow \mathcal{F}(\C^n) $ is a unitary operator which takes the Hermite polynomials $ H_\alpha $ into the monomials $ z^\alpha.$ The inversion formula for $ G $ reads as
$$ G^\ast F(x) = \int_{\C^n} F(w) e^{\frac{1}{4} \bar{w}^2} e^{-i x \cdot \bar{w}} \, d\nu(w).$$
To prove this, let $ F = Gf $ and define $ g(x) = f(x) e^{-\frac{1}{2}|x|^2} $ so that $ F(z) = Bg(iz).$ Consider the Hermite expansion of $ g$:
$$ g(x) = \sum_{\alpha \in \Na^n} (g, \Phi_\alpha) \, \Phi_\alpha(x).$$
Since $ B $ takes $ \Phi_\alpha $ onto $ \zeta_\alpha $ and $ (g, \Phi_\alpha) = (Bg,\zeta_\alpha) $ we get
$$ g(x) = \sum_{\alpha \in \Na^n} (Bg, \zeta_\alpha) \, \Phi_\alpha(x) = \sum_{\alpha \in \Na^n} (F, \zeta_\alpha) \, (-i)^{|\alpha|}\, \Phi_\alpha(x).$$
We now recall that the Hermite polynomials satisfy the generating function identity 
$$ \sum_{\alpha \in \Na^n} \Phi_\alpha(x) \, \zeta_\alpha(w) = \pi^{-n/2} e^{-\frac{1}{2}|x|^2} e^{-\frac{1}{4}w^2+ x \cdot w}.$$
In view of this the  Hermite expansion of $ g $ takes the form
$$ g(x) =  e^{-\frac{1}{2}|x|^2}\, \int_{\C^n} F(w) \,e^{\frac{1}{4}\bar{w}^2} \, e^{-i x \cdot \bar{w}} \, d\nu(w).$$
As $ g(x) = f(x) e^{-\frac{1}{2}|x|^2} $ this proves the inversion formula for the Gauss-Bargmann transform.\\

We record the following properties of the Gauss-Bargmann transform which can be verified directly from the definition using  the relation
$$ 2 \frac{\partial}{\partial z_j} e^{-(x-\frac{i}{2}z)^2} = -i\, \frac{\partial}{\partial x_j} e^{-(x-\frac{i}{2}z)^2}.$$ 

\begin{lem} For any $ f \in L^2(\R^n, d\gamma),$ we have
\begin{enumerate} 
\item $ (-\frac{\partial}{\partial z_j} +\frac{1}{2}z_j)Gf(z) = -i\, G(x_jf)(z) $

\item $ 2 \frac{\partial}{\partial z_j}Gf(z) = i \,G(\frac{\partial}{\partial x_j}f)(z)$

\item $(\frac{\partial}{\partial z_j} +\frac{1}{2}z_j)Gf(z) = -i \,G((-\frac{\partial}{\partial x_j}+x_j)f)(z) $

\end{enumerate}
\end{lem}

\subsection{Hermite, Gauss  and Fock-Sobolev spaces} The explicit spectral decomposition of the Hermite operator $ H = -\Delta+|x|^2 $ on $ \R^n $ allows us to define the fractional powers $ H^s, s \geq 0 $ by
$$  H^sf = \sum_{k=0}^\infty  (2k+n)^s\, P_kf,\,\,\,\,  P_kf = \sum_{|\alpha|=k} (f, \Phi_\alpha)\, \Phi_\alpha.$$
We define the Hermite-Sobolev space $ W_H^{s,2}(\R^n), s\geq 0 $ as the domain of $ H^{s/2}$ and equip it with the inner product $ (f,g)_{s} = (H^{s/2}f, H^{s/2}g).$ When $ s =k $ is a non-negative integer we have
$$ W_H^{k,2}(\R^n) = \{ f \in L^2(\R^n): x^\alpha \partial^\beta f \in L^2(\R^n), |\alpha|+|\beta| \leq k \}.$$
This follows from the well known fact that $ x^\alpha \partial^\beta H^{-(|\alpha|+|\beta|)/2} $ is bounded on $ L^2(\R^n).$ 
The Hermite-Sobolev spaces are invariant under translations and Fourier transform. The latter property is a consequence of the fact that Hermite functions are eigenfunctions of the Fourier transform: $\widehat{\Phi}_\alpha = (-i)^{|\alpha|} \, \Phi_\alpha.$  The Hermite-Sobolev space $ W_H^{s,2}(\R^n) $ is a subspace of the standard Sobolev space $ W^{s,2}(\R^n) $ and hence by Sobolev embedding theorem it follows that $ W_H^{s,2}(\R^n) \subset L^\infty(\R^n) $ whenever $ s > n/2.$ For more about Hermite-Sobolev spaces we refer the reader to \cite{BT}.\\

Let us introduce the annihilation operators $ A_j = \frac{\partial}{\partial x_j}+ x_j $ which take $ \Phi_\alpha $ into $ \sqrt{2\alpha_j } \Phi_{\alpha-e_j}$ and define $ A^\alpha = \Pi_{j=1}^n A_j^{\alpha_j}.$ The formal adjoints $ A_j^\ast = -\frac{\partial}{\partial x_j} + x_j $ are the creation operators which take $ \Phi_\alpha $ into $ \sqrt{2(\alpha_j+1) } \Phi_{\alpha+e_j}.$ Therefore,
$$ W_H^{k,2}(\R^n) = \{ f \in L^2(\R^n): A^\alpha f \in L^2(\R^n), |\alpha| \leq k \}.$$ 
This allows us to identify the Hermite-Sobolev spaces with the Gauss-Sobolev spaces  
$$ W_G^{k,2}(\R^n) = \{ g \in L^2(\R^n, d\gamma):  \partial^\alpha g \in L^2(\R^n, d\gamma), |\alpha| \leq k \}.$$ 
To see this, consider the Hermite (polynomial) coefficients of $ \frac{\partial}{\partial x_j} g \in L^2(\R^n, d\gamma)$: 
$$ \int_{\R^n} \frac{\partial}{\partial x_j}g(x)\, H_\alpha(x)\, d\gamma(x) =  \int_{\R^n} A_j f(x)\, H_\alpha(x)\, e^{-\frac{1}{2}|x|^2}\, dx  $$
where $ f(x) = g(x) e^{-\frac{1}{2}|x|^2}.$ Thus we see that $ \partial^\alpha g \in  L^2(\R^n, d\gamma)$ if and only if $ A^\alpha f \in L^2(\R^n) $ providing a one to one correspondence between these spaces. The same correspondence allows us to define the Gauss-Sobolev spaces for any $ s \geq 0 $ as follows:
$$ W_G^{s,2}(\R^n) = \{ g \in L^2(\R^n, d\gamma):  g(x) e^{-\frac{1}{2}|x|^2} \in  W_H^{s,2}(\R^n) \}.$$ \\

The Fock-Sobolev spaces $ \mathcal{F}^{s,2}(\C^n) $ are defined to be the image of $ W_H^{s,2}(\R^n) $ under the Bargmann transform. Equivalently, we can also define them as the image of 
$W_G^{s,2}(\R^n) $ under the Gauss-Bargmann transform.  It follows from either definition that $ F \in \mathcal{F}^{s,2}(\C^n) $ if and only if 
$$ \sum_{\alpha \in \Na^n} (2|\alpha|+n)^s\, |(F, \zeta_\alpha)|^2 < \infty.$$
In \cite{CP} the authors have shown in  that $ F \in \mathcal{F}^{s,2}(\C^n) $ if and only if 
$$ \int_{\C^n} (1+|z|^2)^s |F(z)|^2\,  d\nu < \infty.$$
From the relation $2 \frac{\partial}{\partial z_j}Gf(z) = i \,G(\frac{\partial}{\partial x_j}f)(z) $ it follows that $ F \in \mathcal{F}^{k,2}(\C^n) $ if and only if $ \partial^\alpha F \in \mathcal{F}(\C^n) $ for all $ |\alpha| \leq k.$ We can also show that this is equivalent to $ z^\alpha F \in \mathcal{F}(\C^n) $ for all $ |\alpha| \leq k.$ 

\subsection{Fourier multipliers on Hermite-Sobolev spaces} In this subsection we gather some information about Fourier multipliers on $ W_H^{s,2}(\R^n).$ Recall that a measurable function $ m $ on $ \R^n $ is said to be a Fourier multiplier on $ W_H^{s,2}(\R^n) $ if the operator $ T_m $ is bounded on that space. We denote by $ MW_H^{s,2}(\R^n) $ the set of all such multipliers. As $ W_H^{s,2}(\R^n) $ is invariant under the Fourier transform, it follows that $ m \in MW_H^{s,2}(\R^n) $ if and only if $ m $ is a pointwise multiplier on the same space. Indeed, this is an immediate consequence of $ \mathcal{F} \circ T_m \circ \mathcal{F}^\ast f = m f.$ This property allows us to get some information on Fourier multipliers for $ W_H^{s,2}(\R^n) .$\\

 It is preferable to get a characterisation of $ MW_H^{s,2}(\R^n) $ at lease when $ s =k $ is an integer. We note that for $ |\alpha| \leq k $ the operators $ f \rightarrow x^\alpha f $ and $ f \rightarrow \partial^\alpha f$ are both bounded from $ W_H^{k,2} \rightarrow W_H^{k-|\alpha|/2}.$ This is a consequence of the fact that   the operators $ x^\alpha H^{-|\alpha|/2} $ and $ \partial^\alpha H^{-|\alpha|/2} $ are bounded on $ L^2(\R^n).$  From the  equation
$$ \partial_j (m f) = (\partial_j m) f + m \partial_j f $$
we infer that $ m \in MW_H^{1,2} $ if and only if $ m$ is bounded and $ (\partial_j m) H^{-1/2} $ is bounded on $ L^2(\R^n).$ Extending this to higher order spaces we get

\begin{lem} For any non-negative integer $ k$, $ m \in MW_H^{k,2}(\R^n) $ if and only if the operators $ (\partial^\alpha m)H^{-|\alpha|/2} $ are bounded on $ L^2(\R^n)$ for all $ |\alpha| \leq k.$
\end{lem}
\begin{proof} For any $ |\mu|+|\alpha| \leq k,$ we have 
$$ x^\mu\, \partial^\alpha (mf)  = \sum_{\beta+\gamma = \alpha} c_{\beta,\gamma} (\partial^\beta m) x^\mu \partial^\gamma f = \sum_{\beta+\gamma =  \alpha} c_{\beta,\gamma} (\partial^\beta m) H^{-|\beta|/2} H^{|\beta|/2} \big(x^\mu \partial^\gamma f) $$
from which the sufficiency follows immediately. For the necessity we use induction on $|\alpha|.$ As noted above the result is true for $ k =1.$ As
$$ c_{\alpha,0}\,( \partial^\alpha m)f  = \partial^\alpha (mf) -\sum_{\beta+\gamma = \alpha, \gamma \neq 0} c_{\beta,\gamma} (\partial^\beta m)  \partial^\gamma f $$
from the assumption that $ m \in MW_H^{k,2}(\R^n), $ the boundedness of $ \partial^\alpha $ and the induction hypothesis, we get the estimate
$$ \| ( \partial^\alpha m)f \|_2  \leq C\, \| f\|_{(k)} = C\,\| H^{|\alpha|/2} f \|_{(k-|\alpha|)} \leq C^\prime \| f\|_2.$$
This proves the necessity part.
\end{proof} 


\section{Fourier Multipliers on Fock spaces} 

\subsection{Boundedness of $ S_\varphi $ and $\widetilde{S}_\varphi$} Let $ S $ be  a bounded linear operator on the Fock space. As observed in the introduction,  with $ K_S(z,\bar{w}) = Sg_{\bar{w}}(z), g_{\bar{w}}(z) = e^{\frac{1}{2}z \cdot \bar{w}},$ we have
$$ SF(z) = \int_{\C^n} F(w)\, K_S(z,\bar{w}) \, d\nu .$$  
If $ S $ commutes with $ \rho(a), a \in \R^n,$ we have 
$$ SF(z+a) e^{-\frac{1}{2} a \cdot z} = \int_{\C^n} F(w+a)\, e^{-\frac{1}{2} a \cdot w}\, K_S(z,\bar{w}) \, d\nu .$$
Therefore, the kernel $K_S $ satisfies the relation 
$$ K_S(z+a,\bar{w}) e^{-\frac{1}{2} a \cdot z} = K_S(z,\bar{w}-a) e^{\frac{1}{2} a\cdot \bar{w}}.$$ For the kernel $ k_S(z,\bar{w}) = K_S(z,\bar{w}) e^{-\frac{1}{2} z\cdot \bar{w}} $ the above relation translates into   $ k_S(z+a, w) = k_S(z,w-a) $ for all $ z, w \in \C^n, a \in \R^n.$ This simply means that $ k_S(z,w) = k_S(z-w,0) $ for all $ z, w \in \R^n$ and the same holds for all $ z, w \in \C^n.$ Therefore, $ K_S(z,\bar{w}) = \varphi(z-\bar{w}) e^{\frac{1}{2} z \cdot \bar{w}} $ for an entire function $ \varphi $ so that
$$ SF(z) = \int_{\C^n} F(w)\, \varphi(z-\bar{w}) e^{\frac{1}{2} z \cdot \bar{w}}\, d\nu .$$
As $e^{\frac{1}{2} z \cdot \bar{w}}$ is the reproducing kernel, it follows that $ \varphi(z) = S1(z) $ and hence $ \varphi \in \mathcal{F}(\C^n).$\\

{\bf Proof of Theorem 1.1} We begin with proving the claimed representation of $ \varphi $ under the assumption that $ S_\varphi $ is bounded.  Since the transferred operator $ T = B^\ast \circ S \circ B $ commutes with translations we have  $ T= T_m $ for a bounded function $ m $ on $ \R^n.$ As $ \varphi(z-\bar{w}) ^{\frac{1}{2}z \cdot \bar{w}} = K_S(z,\bar{w}) = S_\varphi g_{\bar{w}}(z), $ we  will calculate the action of $ S_\varphi  = B \circ T_m \circ B^\ast $ on  $ g_{\bar{w}}(z) = e^{\frac{1}{2} z\cdot \bar{w}} .$  Since the Bargmann transform takes the Hermite functions $ \Phi_\alpha $ onto the monomials $ \zeta_\alpha(z) $ we see that $ g_{\bar{w}}(z) = Bh_{\bar{w}}(z) $ where
$$ h_{\bar{w}}(x) = \sum_{\alpha \in \Na^n} \Phi_\alpha(x) \, \zeta_\alpha(\bar{w}) = \pi^{-n/2} e^{-\frac{1}{2}|x|^2} e^{-\frac{1}{4}\bar{w}^2+ x \cdot \bar{w}}.$$
We also observe that $ \widehat{h}_{\bar{w}}(\xi) = h_{-i\bar{w}}(\xi)$ as $ \Phi_\alpha $ are eigenfunctions of the Fourier transform with eigenvalues $ (-i)^{|\alpha|}.$ Hence we have
$$ B \circ T_m \circ B^\ast g_{\bar{w}}(z) = e^{\frac{1}{4}z^2} \int_{\R^n} m(\xi) e^{-\frac{1}{2}|\xi|^2} e^{\frac{1}{4}\bar{w}^2-i \xi \cdot \bar{w}}  e^{-\frac{1}{2}|\xi|^2} \, e^{i\xi \cdot z} \,d\xi.$$
We can rewrite the above as $ e^{\frac{1}{2} z \cdot \bar{w}} \varphi(z-\bar{w}) $ where
$$ \varphi(z) = \pi^{-n/2} \int_{\R^n} m(\xi) e^{-(\xi-\frac{i}{2} z)^2} d\xi.$$
This proves the statement on the boundedness of $ S_\varphi.$ To deal with $ \widetilde{S}_\varphi $ we proceed as follows.\\

Suppose $ S_\varphi $ corresponds to the Fourier multiplier $ m $ so that $ B \circ T_m \circ B^\ast = S_\varphi .$ Denoting the operator of multiplication by $ m$ simply by $ m $ itself, we 
note that $ \mathcal{F} \circ T_m \circ \mathcal{F}^\ast = m $ and hence $ B \circ m \circ B^\ast =  B \circ \mathcal{F} \circ T_m \circ \mathcal{F}^\ast \circ B^\ast.$  Since $ B \circ \mathcal{F} = U \circ B $ the above gives
$$ B \circ m \circ B^\ast =  U \circ B \circ T_m  \circ B^\ast \circ U^\ast = U \circ S_\varphi \circ U^\ast.$$
A simple calculation shows that
$$ U \circ S_\varphi \circ U^\ast F(z) = \int_{\C^n} F(w) \,e^{\frac{1}{2} z \cdot \bar{w}} \,U\varphi(z+w) \, d\nu(w).$$
Thus we see that the operator $ \widetilde{S}_{U\varphi} $ is bounded on $ \mathcal{F}(\C^n) $ if and only if  $ S_\varphi $ is bounded which happens precisely when
$$ \varphi(z) = \int_{\R^n} m(\xi) e^{-(\xi-\frac{i}{2}z)^2} d\xi .$$ 
Replacing $ \varphi $ by $ U^\ast \varphi$ we  can state this as $ \widetilde{S}_{\varphi} $ is bounded on $ \mathcal{F}(\C^n) $ if and only if  
$$ \varphi(z) = \int_{\R^n} m(\xi) e^{-(\xi-\frac{1}{2}z)^2} d\xi .$$ 
To complete the proof of Theorem 1.1 it remains to be shown that for any non constant $ \varphi $ only one of the operators $ S_\varphi $ or $ \widetilde{S}_\varphi$ can be bounded.
Since the boundedness of $ \widetilde{S}_\varphi $ is equivalent to that of $ S_{U^\ast \varphi}$ as shown above,  we only need to prove  the following proposition which is an easy 
consequence of Hardy's theorem  \cite{H} for the Fourier transform on $ \R^n.$ For the convenience of the readers we recall it here.\\

\begin{thm}  Constant multiples of the Gaussian $ e^{-\frac{1}{2}|x|^2} $ are the only measurable functions $ f $ on $ \R^n $ which satisfy the following estimates:
$$ |f(x)| \leq C\, e^{-\frac{1}{2}|x|^2},\,\,\, |\widehat{f}(\xi)| \leq C\, e^{-\frac{1}{2}|\xi|^2}.$$
\end{thm}  

The above result is an instance of the uncertainty principle for the Fourier transform. For more about Hardy's theorem and it s generalisations to various other setting we refer to the monograph \cite{T2}.\\

\begin{prop} If $ \varphi \in \mathcal{F}(\C^n) $ is non-constant, then  $ S_\varphi $ and $ S_{U^\ast \varphi} $ cannot  both be bounded on $ \mathcal{F}(\C^n).$  Thus  either $ S_\varphi $ or $ \widetilde{S}_\varphi $ is unbounded.
\end{prop}
\begin{proof} Assume the contrary. If $ S_{U^\ast \varphi} $ is bounded then
$$ U^\ast\varphi(z) = e^{\frac{1}{4}z^2} \int_{\R^n} m(\xi) e^{-|\xi|^2} e^{i\, z \cdot \xi} \, d\xi $$ where $ m $ is bounded. This simply means that
$$ \varphi(z) = U^\ast\varphi(-iz) = Bf(z),\,\,  f(x) = m(x) e^{-\frac{1}{2} |x|^2}.$$
Since $ \varphi = U^\ast (U\varphi),$ from the boundedness of $ S_\varphi$ we also have  
$$ U \varphi(z) = Bg(z),\,\,g(x) = \tilde{m}(x) e^{-\frac{1}{2} |x|^2}$$
where $ \tilde{m} $ is also bounded.
Thus $ Bg = U\varphi = B\widehat{f} $ so that $ g = \widehat{f} $ leading to the estimates
$$ |f(x)| \leq C e^{-\frac{1}{2}|x|^2},\,\,\,  |\widehat{f}(\xi)| \leq C e^{-\frac{1}{2}|\xi|^2}.$$
By appealing to Hardy's theorem we conclude that $ f(x) = c\, e^{-\frac{1}{2}|x|^2} $ which is a contradiction since we have assumed that $ \varphi $ is non-constant.\\
\end{proof}

\begin{rem} If we let $ \mathcal{F}_0(\C^n) $  to be  the image of $ L^\infty(\R^n) $ under the Gauss-Bargmann transform, then   $ S_\varphi $ is bounded if and only if $ \varphi \in \mathcal{F}_0(\C^n) .$  In view of this, the  above result can be viewed as an uncertainty principle for the operator $ U $ acting on $ \mathcal{F}(\C^n).$ \begin{it}  For a non constant $ \varphi \in \mathcal{F}(\C^n) $ either  $ \varphi \in \mathcal{F}_0(\C^n) $ or $ U^\ast \varphi  \in \mathcal{F}_0(\C^n) $ but not both.\end{it}\\
\end{rem} 

 The function $m $ which corresponds to $ S_\varphi $ is given by the formula
$$  m(x) =  G^\ast \varphi(x) = \int_{\C^n} \varphi(w) e^{\frac{1}{4} \bar{w}^2} e^{-i x \cdot \bar{w}} \, d\nu(w).$$
A similar formula holds for the function $ m $ associated to $ \widetilde{S}_\varphi.$
For $ \varphi \in \mathcal{F}(\C^n) $ it is clear that $ m \in L^2(\R^n,d\gamma).$ It would be interesting to find conditions on $ \varphi $ so that $ m $ is bounded.

\begin{rem}
Another formula for $ m $ is given in terms of the Taylor coefficients of $ \varphi $:
$$ m(x) e^{-\frac{1}{2}|x|^2}  = \sum_{\alpha \in \Na^n} (\varphi, \zeta_\alpha) \, (-i)^{|\alpha|}\, \Phi_\alpha(x) = 
\sum_{\alpha \in \Na^n}  \frac{(-i)^{|\alpha|}}{ 2^{|\alpha|} \alpha!}  (\varphi, z^\alpha)H_\alpha(x) e^{-\frac{1}{2}|x|^2}, $$
see the proof of the inversion formula for the Gauss-Bargmann transform. From this it is not easy to find conditions  on the coefficients of $ \varphi $ so that $ m $ is bounded. However, for $ m $ to be bounded, $ (\varphi, \zeta_\alpha) $ has to be non zero for infinitely many values of $ \alpha.$
\end{rem}

{\bf Proof of Theorems 1.2 and 1.4:}   In the proof of Theorem 1.1 we have already noted  the relation $ B \circ m \circ B^\ast =  B \circ E \circ m \circ E^\ast \circ B^\ast = U \circ S_\varphi \circ U^\ast ,$ which gives $ G \circ m \circ G^\ast =  S_\varphi  .$  We rewrite this as
$$ G^\ast (S_\varphi F) = m\, G^\ast F = G^\ast \varphi \,  G^\ast F.$$ 
In view of this we can represent the operator $ S_\varphi$ in the  familiar form of a Fourier multiplier:
$$ S_\varphi F(z) = e^{\frac{1}{4}z^2} \, \int_{\R^n} m(x)\, G^\ast F(x) e^{i x \cdot z} \, d\gamma(x) .$$ 
Similarly, the relation $  U \circ S_{U^\ast \varphi} \circ U^\ast = \widetilde{S}_\varphi $ gives us $ B^\ast \circ \widetilde{S}_\varphi \circ B = m $ where $ m $ now is associated to $ U^\ast \varphi,$ 
that is to say $ Gm = U^\ast \varphi $ or $ B (m e^{-\frac{1}{2}|\cdot|^2} ) = \varphi.$ So we can write this as 
$$ B^\ast( \widetilde{S}_\varphi F)(x) = m(x)\, B^\ast F(x) = e^{\frac{1}{2}|x|^2}\, B^\ast \varphi(x)\, B^\ast F(x).$$
This gives the desired representation for the operator $\widetilde{S}_\varphi $:
$$ \widetilde{S}_\varphi F(z) = e^{\frac{1}{4}z^2} \, \int_{\R^n} m(\xi)\, \, B^\ast F(\xi)\, e^{ z \cdot \xi}\, d\xi.$$
In order to complete the proof of Theorem 1.4, we need to show that  any $ \widetilde{S} $ which commutes with $ \rho(ib), b \in \R^n $ has the form $ \widetilde{S}_\varphi $ for some $ \varphi.$
Since $  \mathcal{F} \circ \pi(ia) \circ \mathcal{F}^\ast = \pi(a) $ it follows that $ \widetilde{S} $ commutes with $ \rho(ib) $ if  and only if $ S = U^\ast \circ \widetilde{S} \circ S$ commutes with $ \rho(b) $ and hence the result follows from the argument given at the beginning of this subsection.\\

\subsection{Smoothing properties of $ S_\varphi $ and $ \widetilde{S}_\varphi$}
Let us use the notation $ \varphi \ast F = S_\varphi F .$  For the standard convolution on $ \R^n $ we have
$$ \frac{\partial}{\partial x_j} (g \ast f) = \frac{\partial g}{\partial x_j} \ast f =  g \ast \frac{\partial f}{\partial x_j}.$$
The analogue of this for the above convolution on $ \mathcal{F}(\C^n) $ is given in the following. 

\begin{lem} Let  $ D_j = \frac{\partial}{\partial z_j}+ \frac{1}{2} z_j ,\,  D_j^\ast = -\frac{\partial}{\partial z_j}+ \frac{1}{2} z_j .$ Then for $ \varphi, F \in \mathcal{F}(\C^n),$ we have
\begin{enumerate}
\item $ D_j^\ast (\varphi \ast F) = D_j^\ast  \varphi \ast F = \varphi \ast D_j^\ast F,\,\,$
\item $ D_j (\varphi \ast F) = D_j \varphi \ast F + 2\, \varphi \ast \frac{\partial}{\partial w_j} F = 2  ( \frac{\partial}{\partial z_j} \varphi )\ast F(z) + \varphi \ast D_j F.$
\end{enumerate}
\end{lem}
\begin{proof} We start with

$$ \frac{\partial}{\partial z_j} (\varphi \ast F)(z) = \int_{\C^n} F(w) e^{\frac{1}{2} z \cdot \bar{w}} \big( \frac{1}{2} \bar{w_j}\,   \varphi(z-\bar{w})+ \frac{\partial}{\partial z_j} \varphi(z-\bar{w}) \big)  
d\nu(w).$$
Rewriting  the first integral  as 
$$ - \int_{\C^n} F(w) e^{\frac{1}{2} z \cdot \bar{w}} \,    \varphi(z-\bar{w}) \, \frac{\partial}{\partial w_j} e^{- \frac{1}{2}\sum_{j=1}^n w_j \bar{w}_j}\,  dw $$
and integrating by parts we obtain 
$$ \frac{\partial}{\partial z_j} (\varphi  \ast F) =   \frac{\partial}{\partial z_j} \varphi  \ast F+  \varphi \ast \frac{\partial}{\partial w_j} F.$$
We can also prove the relation 
$$ \frac{1}{2} z_j (\varphi \ast F)(z) =   (\frac{1}{2} z_j \,\varphi) \ast  F(z) +  \varphi \ast  \frac{\partial}{\partial w_j} F(z)$$
by starting with the following and proceeding as above:
$$ \frac{1}{2} z_j S_\varphi F(z) = \int_{\C^n} F(w) e^{\frac{1}{2} z \cdot \bar{w}} \big(\frac{1}{2} (z_j-\bar{w}_j) \varphi(z-\bar{w}) +  \frac{1}{2} \bar{w_j}\,   \varphi(z-\bar{w}) \big)  d\nu(w).$$ Integrating by parts in the integral
$$ S(\varphi, \frac{1}{2}w_j F) =  - \int_{\C^n} F(w) e^{\frac{1}{2} z \cdot \bar{w}} \,    \varphi(z-\bar{w}) \, \frac{\partial}{\partial \bar{w}_j} e^{- \frac{1}{2}\sum_{j=1}^n w_j \bar{w}_j}\,  dw $$
we also obtain
$$ \frac{1}{2} z_j (\varphi \ast F)(z) =  ( \frac{\partial}{\partial z_j} \varphi )\ast F(z) +   \varphi \ast (\frac{1}{2}w_j F)(z).$$
The lemma is proved by combining these three relations.\\
\end{proof}

The behaviour the convolution on $ \mathcal{F}(\C^n) $ under the derivatives established in the above lemma has some interesting consequences on the smoothing properties of the operator 
$ S_\varphi.$  Let us define the subspace $ \mathcal{F}_{1,0}^{k,2}(\C^n)$ by 
$$ \mathcal{F}_{1,0}^{k,2}(\C^n) = \{ F \in \mathcal{F}(\C^n) : D^{\ast \alpha} F \in \mathcal{F}(\C^n),\, |\alpha|\leq k  \}.$$
We define $ \mathcal{F}_{0,1}^{k,2}(\C^n)$ in a similar fashion using the operators $ D_j.$ It follows from the definition of the Bargmann transform that
$ D_j Bf(z) = B(\xi_j f)(z).$ From the equivalent definition of the Bargmann transform, namely
$$ Bf(z) =  e^{\frac{1}{4}z^2} \int_{\R^n} \widehat{f}(\xi) e^{-\frac{1}{2}|\xi|^2} e^{i \xi \cdot z} d\xi,$$
we also have   $ D_j^\ast Bf(z) = B (\partial_j f )(z) .$ Hence $  \mathcal{F}_{1,0}^{k,2}(\C^n)$ is the image of the standard Sobolev space $ W^{k,2}(\R^n) $ under the Bargmann transform
 whereas $  \mathcal{F}_{0,1}^{k,2}(\C^n)$ is the image of $ \mathcal{F}W^{k,2}(\R^n), $ the Fourier image of the Sobolev space $ W^{k,2}(\R^n)$ under the Bargmann transform. Thus $ f \in W^{k,2}(\R^n)$ if and only if $ \widehat{f} \in  \mathcal{F}W^{k,2}(\R^n), $ and hence $ F \in  \mathcal{F}_{1,0}^{k,2}(\C^n)$ if and only if $ UF \in   \mathcal{F}_{0,1}^{k,2}(\C^n).$
From Lemma 3.1 we immediately get

\begin{prop} For any $ \varphi \in \mathcal{F}_0(\C^n),$ the operator $ S_\varphi $ is bounded on $  \mathcal{F}_{1,0}^{k,2}(\C^n). $  If we assume that $ D^{\ast \alpha}\varphi \in \mathcal{F}_0(\C^n)$ for all $ |\alpha| \leq k,$ then $ S_\varphi : \mathcal{F}(\C^n) \rightarrow   \mathcal{F}_{1,0}^{k,2}(\C^n) $  is bounded. Under the same assumption on $ \varphi $ the operator $ U \circ S_{\varphi} \circ U^\ast = \widetilde{S}_{ U\varphi} : \mathcal{F}(\C^n) \rightarrow  \mathcal{F}_{0,1}^{k,2}(\C^n)$ is also bounded.\\
\end{prop} 

\begin{rem} If $ m $ is the multiplier associated to $ S_\varphi,$ then from Lemma 2.1 we have $ D_j^\ast \varphi = -i\, G(x_j m) $ and hence we can state the above result as follows: if we assume that $ x^\alpha m(x) $ is bounded for all $ |\alpha| \leq k,$ then $ S_\varphi : \mathcal{F}(\C^n) \rightarrow   \mathcal{F}_{1,0}^{k,2}(\C^n) $  is bounded.\\
\end{rem}

\begin{prop} Suppose $ \partial^\alpha \varphi \in \mathcal{F}_0(\C^n) $ for all $ |\alpha| \leq k.$ Then $ S_\varphi $  is bounded on $  \mathcal{F}_{1,0}^{k,2}(\C^n). $ If the multiplier $ m $  associated to $ S_\phi$ satisfies the condition $ \partial^\alpha m \in L^\infty(\R^n) $ for all $ |\alpha| \leq k,$ then $ S_\varphi $  is bounded on $  \mathcal{F}_{1,0}^{k,2}(\C^n). $
\end{prop}
\begin{proof} The first statement above  is an immediate consequence of  the following relation which follows from (2) of Lemma 3.1:
$$ D^\alpha (\varphi \ast F) =  \sum_{\beta + \gamma = \alpha} c_{\beta,\gamma} \,   \partial^\beta \varphi  \ast D^\gamma F.$$ 
The condition on $ m = G^\ast \varphi $ can be translated into the condition on $ \varphi$ using 
$$ 2 \frac{\partial}{\partial z_j}Gm(z) = i \,G(\frac{\partial}{\partial x_j}m)(z)$$
proved in Lemma 2.1. This proves the second statement.\\
\end{proof}

Unlike the spaces $\mathcal{F}_{1,0}^{k,2}(\C^n) $ and $\mathcal{F}_{0,1}^{k,2}(\C^n),$ the spaces $\mathcal{F}^{s,2}(\C^n) $ are invariant under the action of $ U$ and hence  we have the results stated in Theorems 1.5 and 1.7.\\

{\bf Proofs of Theorems 1.5 and 1.7:}  Theorem 1.5 follows from the fact that $ S_\varphi$ is bounded on $\mathcal{F}^{s,2}(\C^n) $ if and only if $ T_m = B^\ast \circ S_\varphi \circ B $ is bounded on $ W_H^{s,2}(\R^n).$  When $ s = k$ is a non-negative integer and $ \varphi \in \mathcal{F}_0^{k,2}(\C^n)$ it follows from the definition that $ m = G^\ast \varphi \in L_k^\infty(\R^n) $ which defines a multiplier on $ W_H^{k,2}(\R^n).$

\begin{rem} As noted earlier $ \partial^\alpha \varphi = G(\partial^\alpha m).$ It would be interesting to find conditions on the derivatives $ \partial^\alpha \varphi $ so that $ m \in L_k^\infty(\R^n).$ We do not know how to do this.
\end{rem}

\section{Pseudo-differential operators on the Fock space} In this section we prove Theorem 1.11 and and calculate the kernels $ K $ associated to a class of radial symbols $ \sigma.$ 

\subsection{Proof of Theorem 1.11}
 We  show that  if $ K(z,w) = B\psi(z,w),$ then the transferred operator $ B \circ W(\sigma)  \circ B^\ast  = S_K $  where $ \sigma $ is given by \eqref{sigma}.
Consider the action of  $ B \circ W(\sigma) \circ B^\ast $ on the function $ g_w(z) = e^{\frac{1}{2} z \cdot w} = B h_w(z)$. Recalling the definition of $ B $ and $ W(\sigma) $ we get,
$$ B\circ W(\sigma)\circ B^\ast g_w(z) = e^{\frac{1}{4}z^2} \int_{\R^n}\Big( \int_{\R^{2n}} \sigma(x,y) e^{i(x\cdot \xi+\frac{1}{2} x\cdot y)} h_w(\xi+y)  dx dy\Big) e^{-\frac{1}{2}(\xi-z)^2} d\xi.$$
The explicit formula for $ h_w(x) $ and the change of variables $ \xi \rightarrow \xi-y/2$ simplifies the above to give
$$e^{\frac{1}{4}(z^2+w^2)} \int_{\R^{2n}}\Big( \int_{\R^{n}} \sigma(x,y) e^{i x\cdot \xi}   dx \Big) e^{-\frac{1}{2}(\xi+y/2-w)^2} e^{-\frac{1}{2}(\xi-y/2-z)^2} dy d\xi.$$
Another change of variables $ u = \xi+y/2, v = \xi-y/2 $ reduces the above integral to
$$e^{\frac{1}{4}(z^2+w^2)} \int_{\R^{2n}} \mathcal{F}_1^{-1}\sigma((u+v)/2, u-v) e^{-\frac{1}{2}(u-w)^2} e^{-\frac{1}{2}(v-z)^2} dy d\xi.$$
This proves that $B\circ W(\sigma)\circ B^\ast g_w(z) = B\psi(w,z) $ where $ \psi(u,v) = \mathcal{F}_1^{-1}\sigma((u+v)/2, u-v) .$ As the operator is uniquely determined by  its action on $ g_w $ we get the result.\\

\begin{rem} We now consider the boundedness of the operators $ S_K = B \circ W(\sigma) \circ B^\ast $ on the Fock-Sobolev spaces which is equivalent to the boundedness of $ H^{s/2} \circ W(\sigma) \circ H^{-s/2} $ on $ L^2(\R^n).$ The operator $ H^{s/2} $ is a pseudo-differential operator whose Weyl symbol belongs to the Shubin class $ \Gamma^s(\R^{2n}).$  We say that $ a_s $ is the Weyl symbol of $ H^{s/2} $ if  $ W(\widehat{a}_s) = H^{s/2}.$ The Shubin class $ \Gamma^s(\R^{2n}) $ is defined be the set  of all $ a \in C^\infty(\R^{2n}) $ satisfying the uniform estimates
$$ |\partial_x^\alpha \partial_y^\beta a(x,y)| \leq C_{\alpha,\beta} (1+|x|+|y|)^{s-|\alpha|-|\beta|},\,\,\, \alpha, \beta \in \Na^n.$$
It is known that the Weyl symbol $ a_s $ of $ H^{s/2} $ belongs to $ \Gamma^s(\R^{2n}),$ see Proposition 2.3 in \cite{BMNTT} and Theorem 2.1 in \cite{T3}.
If we assume that $ \sigma $ is such that $ \partial_x^\alpha \partial_y^\beta \widehat{\sigma} $ are bounded for all $ \alpha $ and $\beta,$ then $B \circ W(\sigma) \circ B^\ast = W(a) $
where the Weyl symbol $ \widehat{a} $ also satisfies  $ | \partial_x^\alpha \partial_y^\beta \widehat{a}(x,y)| \leq C_{\alpha,\beta}.$  Therefore, by Calderon-Vaillancourt theorem we can conclude that 
$ B \circ W(\sigma) \circ B^\ast $ is bounded on $ L^2(\R^n).$ For more about pseudo-differential operators we refer to the book by Nicola and Rodino \cite{NR}.
\end{rem}

\subsection{Pseudo-differential operators with radial symbol}  For all the results concerning Laguerre expansions and Weyl transform stated in this subsection we refer to \cite{T1} and \cite{T2}. Let  $ L_k^{n-1}(t) $ be the Laguerre polynomials of type $(n-1)$ and define
$$ \varphi_k^{n-1}(z) = L_k^{n-1}(\frac{1}{2}|z|^2) e^{-\frac{1}{4}|z|^2}.$$ Then it is known that any radial function $ \sigma \in L^1(\C^n)$   has a Laguerre expansion
$$ \sigma(z) = (2\pi)^{-n} \sum_{k=0}^\infty R_k^{n-1}(\sigma)\, \varphi_k^{n-1}(z).$$
Here $ R_k^{n-1}(\sigma) $ are the Laguerre coefficients of $ \sigma$ defined by
$$ R_k^{n-1}(\sigma) = \frac{k! (n-1)!}{(k+n-1)!} \int_{\C^n} \sigma(z) \varphi_k^{n-1}(z)\, dz .$$ 
Moreover, it is also known that $ W(\varphi_k^{n-1}) = (2\pi)^n\, P_k $ where $ P_k $ are the Hermite projections taking $ L^2(\R^n) $ onto the span of $ \Phi_\alpha, |\alpha| =k.$ Consequently, 
$$ W(\sigma) = \sum_{k=0}^\infty R_k^{n-1}(\sigma)\,P_k $$ 
which is a bounded linear operator on $ L^2(\R^n) $ since $ |R_k^{n-1}(\sigma)| \leq c\, \| \sigma\|_1 $ which follows from the fact that $ |\varphi_k^{n-1}(z)| \leq \, \frac{(k+n-1)!}{k!(n-1)!}.$ The same conclusion holds under the weaker assumption that $ \sigma $ is a radial Borel measure on $ \C^n $ for which that coefficients $ R_k^{n-1}(\sigma) $ are uniformly bounded.\\

When $ \sigma $ is as above, the kernel of the operator $ B \circ W(\sigma) \circ B^\ast $ can be calculated explicitly in terms of $ \sigma.$ We have

\begin{prop} Let $ \sigma $ be a radial symbol for which $ W(\sigma) $ is bounded. Then $ B \circ W(\sigma) \circ B^\ast = S_K $ where the kernel is given by the formula
$$ K(z,w) = c_n \,e^{\frac{1}{2} z\cdot w} \, \int_{\C^n}  \sigma(\zeta)\, \frac{J_{n-1}(|\zeta| \sqrt{ z \cdot w})}{(|\zeta| \sqrt{z \cdot w})^{n-1}}\, e^{-\frac{1}{4}|\zeta|^2} \, d\zeta.$$
\end{prop}
\begin{proof} As before $ K(z,w) = (B \circ W(\sigma) \circ B^\ast) g_w(z) = (B \circ W(\sigma)) h_w(z).$ The expansion of $ h_w$ in terms of Hermite functions and the formula for $ W(\sigma) $ together give us
$$ W(\sigma)h_w(x) = \sum_{k=0}^\infty R_k^{n-1}(\sigma)\, \big( \sum_{|\alpha|=k} \Phi_\alpha(x)\,\zeta_{\alpha}(w)\big).$$
Applying the operator $ B $ to the function $ W(\sigma)h_w $ and using the fact that $ B \Phi_\alpha = \zeta_\alpha$ we get
$$ K(z,w) = \pi^{-n/2}\, \sum_{k=0}^\infty R_k^{n-1}(\sigma)\, \big( \sum_{|\alpha|=k} \frac{ z^\alpha\,w^\alpha}{2^{|\alpha|} \alpha !}\big)  = \pi^{-n/2}\, \sum_{k=0}^\infty  \, R_k^{n-1}(\sigma)\, \frac{ ( z\cdot w)^k}{2^k k!}.$$
Recalling the definition of $ R_k^{n-1}(\sigma) $ we have proved that
$$ K(z,w) = c_n \int_{\C^n} \Big( \sum_{k=0}^\infty \frac{L_k^{n-1}(\frac{1}{2}|\zeta|^2)}{(k+n-1)!}\, (\frac{1}{2} z\cdot w)^k  \Big)\, e^{-\frac{1}{4}|\zeta|^2} \sigma(\zeta) \,d\zeta. $$
The proof is completed by using the following generating function identity satisfied by Laguerre polynomials $ L_k^\delta(r) $ of type $ \delta >-1$: for any  $ r >0, t \in \C $ 
$$ \sum_{k=0}^\infty \frac{L_k^{\delta}(r)}{\Gamma(k+\delta+1)!}\, t^k  = e^t\, (\sqrt{rt})^{-\delta}\, J_\delta(2 \sqrt{rt}) $$
where $ J_\delta(t) $ is the Bessel function of order $ \delta,$ see \cite{Sz}.
\end{proof}

\begin{rem} If we write $ \sigma(\zeta) = \sigma_0(|\zeta|) $ then the kernel of the operator $ B \circ W(\sigma) \circ B^\ast $ is given by $ K(z,w) = \mathcal{B}_{n-1} \sigma_0(\sqrt{z\cdot w})$
where for any $ \delta \geq -1/2$ and $ g \in L^2(\R^+, r^{2\delta+1} dr),$
$$ \mathcal{B}_{\delta}g(\tau) = e^{\frac{1}{2}\tau^2} \int_0^\infty  g(r)\, e^{-\frac{1}{4}r^2}\, \frac{J_{\delta}(r\tau)}{(r\tau)^{\delta}} \, r^{2\delta+1} dr.$$
These operators $ \mathcal{B}_\delta,$ known  as generalised Bargmann transforms have been studied in the literature, see \cite{C}. The image of $L^2(\R^+, r^{2\delta+1} dr) $ under $ \mathcal{B}_\delta $ can be characterised as a Hilbert space of even entire functions on $ \C$ which are square integrable with respect to a weight function $ w_\delta(\tau),$ which is known explicitly, expressible in terms of the Macdonald function $K_\delta.$
\end{rem}

\begin{center} \bf{Acknowledgments}
\end{center}
This work was carried out while the author was visiting the University of Queensland, Brisbane  during December 2022- January 2023. He wishes to thank the Department of Mathematics for the hospitality and  Prof. Dietmar Oelz for the invitation.

\end{document}